\documentclass{amsart}
\usepackage{amssymb}

\theoremstyle{plain}

\numberwithin{equation}{section}

\newtheorem*{teorema}{Theorem A}

\newtheorem{lem}[equation]{Lemma}

\newtheorem{claim}[equation]{Claim}

\newcommand{\Irr}{\operatorname{Irr}}
\newcommand{\Lin}{\operatorname{Lin}}
\newcommand{\Z}{\operatorname{Z}}

\newcommand{\modu}{\operatorname{mod}}

\newcommand{\core}{\operatorname{core}}
\newcommand{\Ker}{\operatorname{Ker}}

\theoremstyle{definition}
\newtheorem{comment}[equation]{}

\begin{document}	
\title{Characters of prime degree}

\author{Edith Adan-Bante}

\address{Department of Mathematical Science, Northern Illinois University, 
 Watson Hall 320
DeKalb, IL 60115-2888, USA}

\email{EdithAdan@illinoisalumni.org}


\subjclass{20c15}

\date{2008}

\begin{abstract}
Let $G$ be a finite nilpotent group, $\chi$ and $\psi$ be irreducible complex characters of $G$ of 
prime degree. Assume that  $\chi(1)=p$. Then either the product $\chi\psi$ is a multiple of an irreducible 
character or $\chi\psi$ is the linear combination of at least $\frac{p+1}{2}$ distinct 
irreducible characters. 
\end{abstract}

\maketitle

\begin{section}{Introduction}
Let $G$ be a finite group and $\chi,\psi\in \Irr(G)$ be irreducible complex characters of $G$. 
We can check that the product of $\chi\psi$, where $\chi\psi(g)=\chi(g)\psi(g)$ for all 
$g$ in $G$, is a character and so it can be expressed as a 
linear combination of irreducible characters. Let $\eta(\chi\psi)$ be the number
of distinct irreducible constituents of the product $\chi\psi$.
 
 Fix a prime number $p$. Assume that $G$ is a nilpotent group and
 $\chi$, $\psi\in \Irr(G)$ are characters of degree $p$. In this note, we study
 the product $\chi\psi$.  
 
\begin{teorema}
  Let $G$ be a finite nilpotent group, 
  $\chi$ and $\psi$ be irreducible complex characters of prime degree. Assume that $\chi(1)=p$. Then one
  of the following holds:

  (i) $\chi\psi$ is the sum of $p^2$ distinct linear characters.

  (ii) $\chi\psi$ is the sum of $p$ distinct linear characters of $G$ and of $p-1$
  distinct irreducible characters of $G$ with degree $p$.

  (iii) all the irreducible constituents of $\chi\psi$ are of degree $p$. Also, either
  $\chi\psi$ is a multiple of an irreducible character, or it has at least $\frac{p+1}{2}$
  distinct irreducible constituents, and at most $p$ distinct irreducible
  constituents, i.e. \begin{equation*}
  \mbox{ either }\eta(\chi\psi)=1 \mbox{ or } \frac{p+1}{2}\leq \eta(\chi\psi)\leq p.
  \end{equation*}
  
   (iv) $\chi\psi$ is an irreducible character.
  \end{teorema}
  
  It is proved in Theorem A of \cite{edith1} that given any prime $p$, any $p$-group $P$, any faithful
  characters $\chi,\psi\in \Irr(P)$, either the product
 $\chi\psi$ is a multiple of
 an irreducible, or $\chi\psi$ is the linear combination of at least $\frac{p+1}{2}$ distinct irreducible 
 characters, i.e. either $\eta(\chi\psi)=1$ or $\eta(\chi\psi)\geq \frac{p+1}{2}$. 
 It is proved in \cite{restriction} 
 that given any prime $p$ and any integer $n>0$, there exists a $p$-group $P$ and characters 
  $\varphi,\gamma\in \Irr(P)$ such that $\eta(\varphi\gamma)=n$. Thus without the hypothesis that
  the characters in Theorem A of \cite{edith1} are faithful, the result may not hold.
  In this note, we are proving
  that if the characters have ``small" degree  then the values that
  $\eta(\chi\psi)$ can take have the same constraint as if they were faithful.
  \end{section}
  \begin{section}{Proofs}
  We are going to use the notation of \cite{isaacs}. In addition, we 
  denote by $\Lin(G)=\{\chi\in \Irr(G)\mid \chi(1)=1\}$ the set of linear characters,
  and by $\Irr(G\modu N)=\{\chi\in \Irr(G)\mid \Ker(\chi)\geq N\}$ the set of irreducible
  characters of $G$ that contain in their kernel the subgroup $N$. 
  Also, denote by $\overline{\chi}$ the complex conjugate of $\chi$, i.e.
  $\overline{\chi}(g)=\overline{\chi(g)}$ for all $g$ in $G$.
  
  \begin{lem}\label{derived2}
 Let $G$ be a finite group and $\chi,\psi \in \Irr(G)$. Let $\alpha_1,\alpha_2,
 \ldots, \alpha_n$, for some $n>0$, be the distinct irreducible constituents 
 of the product $\chi\psi$ and $a_1,a_2,\ldots, a_n$ be the unique positive
 integers such that
\begin{equation*}\label{sumapsi}	
 \chi\psi=  \sum_{i=1}^n a_i \alpha_i.
\end{equation*}
 
If $\alpha_1(1)=1$, then $\psi\overline{\alpha_1}=\overline{\chi}$.
Hence  the distinct irreducible constituents of 
the character $\chi\overline{\chi}$ are 
$1_G$, $ \overline{\alpha_1}\alpha_2$, $\overline{\alpha_1}\alpha_2, \ldots, 
\overline{\alpha_1}\alpha_n$, and 
\begin{equation*}
 \chi\overline{\chi}= a_1 1_G + \sum_{i=2}^n a_i (\overline{\alpha_1}\alpha_i).
\end{equation*}
\end{lem}
  \begin{proof} See Lemma 4.1 of \cite{edithderived2}.
  \end{proof}
  
  \begin{lem}\label{linearconstituent}
Let $G$ be a finite $p$-group for some prime $p$ and 
$\chi\in \Irr(G)$ be a character of degree $p$.
Then one of the following holds:

  (i) $\chi\overline{\chi}$ is the sum of $p^2$ distinct linear characters.

  (ii) $\chi\overline{\chi}$ is the sum of $p$ distinct linear characters of $G$ and of $p-1$
  distinct irreducible characters of $G$ with degree $p$.
  \end{lem}
  \begin{proof}
  See Lemma 5.1 of \cite{edith2}.
  \end{proof}

  \begin{lem}\label{basicolemma}
 Let $G$ be a finite $p$-group, for some prime $p$, and $\chi, \psi\in \Irr(G)$ be characters of 
 degree $p$. Then either $\eta(\chi\psi)=1$ or $\eta(\chi\psi)\geq \frac{p+1}{2}$.
  \end{lem} 
\begin{proof}
Assume that the lemma is false. Let $G$ and $\chi,\psi\in \Irr(G)$ be a counterexample of the statement, i.e. $\chi(1)= \psi(1)=p$ and $1<\eta(\chi\psi)<\frac{p+1}{2}$. 

Working with the group $G/(\Ker(\chi)\cap \Ker(\psi))$,
by induction on the order of $G$, we may assume that $\Ker(\chi)\cap \Ker(\psi)=\{1\}$.
Set $n=\eta(\chi\psi)$.
Let $\theta_i\in \Irr(G)$, for  $i=1,\ldots, n $, be the  
 distinct
irreducible
constituents of $\chi \psi$. Set  
\begin{equation}\label{product}
\chi \psi = \sum_{i=1}^n m_i \theta_i
\end{equation}
\noindent where $m_i>0$ is the multiplicity of $\theta_i$ in $\chi\psi$.

If $\chi\psi$ has a linear constituent, then by Lemma \ref{derived2}
and Lemma \ref{linearconstituent} we have that $\eta(\chi\psi)\geq p$.
If $\chi\psi$ has an irreducible constituent of degree $p^2$, then $\chi\psi\in \Irr(G)$
and so $\eta(\chi\psi)=1$.
Thus we may assume that $\theta_i(1)=p$ for $i=1,\ldots, n$.

Since $G$ is a $p$-group, there must exist a subgroup 
$H$ and a linear character $\xi$ of $H$ such that $\xi^G=\chi$. Then
$|G:H|=\chi(1)=p$ and thus $H$ is a normal subgroup. By Clifford theory we have then
\begin{equation}\label{chiintermsofxi}
\chi_H=\sum_{i=1}^p \xi_i
\end{equation}
\noindent for some $\xi_1=\xi, \ldots, \xi_p$ distinct linear characters of $H$.

\begin{claim}\label{inducing}
$H$ is an abelian group.
\end{claim}
\begin{proof}
Suppose that $\psi_H\in\Irr(H)$. Since $(\xi\psi_H)^G=\chi\psi$
by Exercise 5.3 of \cite{isaacs},
 and $\xi\psi_H
\in \Irr(H)$, it follows that either $\xi\psi_H$ induces irreducibly, and thus
$\eta(\chi\psi)=1$, or $\xi\psi_H$ extends to $G$ and thus $(\xi\psi_H)^G$ is the 
sum of the $p$ distinct extensions of $\xi\psi_H$, i.e. 
$\eta(\chi\psi)=p$. Therefore $\psi_H\not\in \Irr(G)$ and since $H$ is 
normal in $G$ of index $p$ and $\psi(1)=p$, $\psi$ is induced from some $\tau\in\Lin(H)$.

Since $H$ is normal of index $p$, both $\xi$ and $\tau$ are linear characters and thus
$\Ker(\xi)\cap \Ker(\tau)\geq [H,H]$. Observe that 
$\core_G(\Ker(\xi)\cap \Ker(\tau))= \core_G(\Ker(\xi))\cap \core_G(\Ker(\tau))=
\Ker(\chi)\cap \Ker(\psi)$. Since $H$ is a normal subgroup of $G$, so is $[H,H]$ and thus
$\{1\}=\Ker(\chi)\cap \Ker(\psi)\geq [H,H]$. Therefore 
$H$ is abelian.
 \end{proof}
 
By the previous claim, observe that also $\psi$ is induced by some linear character $\tau$ of 
$H$ and thus
\begin{equation}\label{psiintermsoftau}
\psi_H=\sum_{i=1}^p \tau_i
\end{equation}
\noindent for some $\tau_1=\tau, \ldots, \tau_p$ distinct linear characters of $H$. 
Observe also that 
 the center of both $\chi$ and $\psi$ is contained in $H$. 

\begin{claim} ${\bf Z}(G)={\bf Z}(\chi)={\bf Z}(\psi)$.
\end{claim}
\begin{proof} 
Suppose that ${\bf Z}(\chi)\neq {\bf Z}(\psi)$.
Set $U={\bf Z}(\chi)\cap {\bf Z}(\psi)$. Either $U$ is properly contained in ${\bf Z}(\chi)$, or
it is properly contained in ${\bf Z}(\psi)$. We may assume that $U< {\bf Z}(\psi)$ and thus
we may find a subgroup 
 $T\leq {\bf Z}(\psi)$ such that $T/U$ is 
chief factor of $G$.  Since $H$ is abelian, ${\bf Z}(\psi)<H$ and $\tau^G=\psi$, then $\psi_T=p\tau_T$ and so $(\tau_i)_T=\tau_T$ for $i=1,\ldots, p$.
 Because $\xi^G=\chi$, $\xi\in \Lin(H)$ and 
$T\not\leq {\bf Z}(\chi)$, 
the stabilizer of $\xi_T$ is $H$. Thus the stabilizer of $\xi_T\tau_T$ in $G$ is $H$.
 By Clifford theory we have that $\xi\tau_i\in \Lin(H)$ induces
irreducibly and $\xi\tau_i$ are distinct characters for $i=1,\ldots, p$. By \eqref{psiintermsoftau} we have that
$\chi\psi=(\xi \psi_H)^G= (\xi(\tau_1+ \cdots+\tau_p))^G=(\xi\tau_1)^G +\cdots+(\xi\tau_p)^G$, and thus
$\eta(\chi\psi)=p$. We conclude that such $T$ can not exist and so 
$ {\bf Z}(\chi)={\bf Z}(\psi)$.

Given any $z \in \Z(\chi)$ and $g \in G$, we have $z^g \cong z \pmod{\Ker(\chi)}$ since $\Z(G/\Ker(\chi)) = \Z(\chi)/\Ker(\chi)$. Hence $[z, g] = z^{-1}z^g$ lies in $\Ker(\chi)$. This same $z$ lies in $\Z(\psi) = \Z(\chi)$. Hence $[z,g]$ also lies in $\Ker(\psi)$. Therefore $[z,g] \in \Ker(\chi) \cap \Ker(\psi) = 1$ for every $z \in \Z(\chi) = \Z(\psi)$ and every $g \in G$. This implies that $\Z(\chi) = \Z(\psi) = \Z(G)$.
\end{proof}

Set $Z={\bf Z}(G)$.  Since $Z$ is the center of $G$,
$\xi^G=\chi$ and $\tau^G=\psi$, we have
\begin{equation}\label{restriction}
\chi_Z= p \xi_Z \mbox{ and } \psi_Z= p \tau_Z.
\end{equation}
Because $\chi_Z\psi_Z=p^2 \xi_Z\tau_Z$,  \eqref{product} implies that 
\begin{equation}\label{thetai}
(\theta_i)_Z = p\xi_Z \tau_Z
\end{equation}
\noindent for all $i=1,\ldots,n$.

Let $Y/Z$ be a chief factor of $G$ with $Y\leq H$.
Since $Z$ is the center of $G$ and $Z={\bf Z}(\chi)$, 
the set $\Lin(\, Y \mid \xi_Z\,)$ of all extensions 
of $\xi_Z$ to linear characters is  $\{(\xi_1)_Y=\xi_Y,(\xi_2)_Y, \ldots, (\xi_p)_Y\}$ and it is a single
 $G$-conjugacy class. By Clifford theory
we have that 
\begin{equation}\label{restrictionchi}
\chi_Y=  \sum_{i=1}^p (\xi_i)_Y.
\end{equation}
Since $H$ is the stabilizer of $\tau_Y$ in $G$ and $\psi(1)=p$,  as before we have that
the set $\Lin(\, Y\mid \tau_z\,)=\{(\tau_1)_Y=\tau_Y,(\tau_2)_Y, \ldots, (\tau_p)_Y\}$ and 
\begin{equation}\label{restrictionpsi}
\psi_Y=  \sum_{i=1}^p (\tau_i)_Y.
\end{equation}

\begin{claim}\label{faithful}
The stabilizer $\{ g \in G \mid (\xi_Y\tau_Y)^g=\xi_Y\tau_Y \}$ 
 of $\xi_Y\tau_Y\in \Lin(Y)$ in $G$ is $H$. 
\end{claim}
\begin{proof}
Assume notation \eqref{product}.
Suppose $\xi_Y\tau_Y$ is a
$G$-invariant character. 
Since $|Y:Z|=p$ and $\xi_Y\tau_Y$ is an extension of $\xi_Z\tau_Z$,
it follows then that 
all the extensions of $\xi_Z\tau_Z$ to $Y$ are $G$-invariant. Thus
by \eqref{product} and \eqref{thetai}, given any $i$, there exists some extension
$\upsilon_i\in \Lin(Y)$ of $\xi_Z\tau_Z$ such that 
 $(\theta_i)_Y=p\upsilon_i$.
  Thus $(\chi\psi)_Y= (\sum_{i=1}^n m_i\theta_i)_Y= \sum_{i=1}^n m_i (\upsilon_i)_Y= \sum_{i=1}^n m_i p \upsilon_i$ has 
at most $n<\frac{p+1}{2}$ distinct irreducible constituents. On the other hand, by
\eqref{restrictionchi} and \eqref{restrictionpsi} we have
\begin{equation*}
(\chi\psi)_Y = \chi_Y\psi_Y= 
(\sum_{i=1}^p (\xi_i)_Y) (\sum_{j=1}^p (\tau_j)_Y) = p
\sum_{j=1}^p \xi_Y(\tau_j)_Y,
\end{equation*}
\noindent and so $(\chi\psi)_Y$ has $p$ distinct irreducible constituents. That is a contradiction and thus $G_{\xi_Y\tau_Y}=H$.
 \end{proof}

By Clifford theory and the previous claim, we have that for each $i=1, \ldots, n$,
there exists a  unique character
 $\sigma_i \in \Lin(\,H \mid  \xi_Y\tau_Y\,)$ such that
 \begin{equation}\label{theta_i}
\theta_i=(\sigma_i)^G.
\end{equation}

If $Y=H$ then $|G: Z|=|G:H||H:Z|= p^2$. Since
$\chi(1)=\psi(1)=p$, by Corollary 2.30 of \cite{isaacs}
 we have that $\chi$ and $\psi$ vanish outside
$Z$. Since $\theta_i(1)=p$ for all $i$ and $|G:Z|=|G:{\bf Z}(\theta_i)|=p^2$, it follows that there exists a unique irreducible character lying above $\xi_Z\tau_Z$
and thus $\eta(\chi\psi)=1$.

\begin{comment}\label{xcentralh}
Fix a subgroup  $X\leq H$ of $G$ such that $X/Y$ is a chief factor
of $G$. Let
 $\alpha$, $\beta \in \Lin (X)$ 
be the  linear characters  such that
\begin{equation*}
\xi_X= \alpha \mbox{ and }
 \tau_X=  \beta.
\end{equation*}  
Since $\sigma_i$ lies above $\xi_Y\tau_Y\in \Lin(Y)$  for
all $i$ and 
$X/Y$ is a chief factor of a $p$-group,  there is
some $\delta_i \in \Irr( X \mod Y)$ such that
\begin{equation}\label{thetaidelta}
(\sigma_i)_X = \delta_i\alpha\beta.
\end{equation}
\end{comment} 
\begin{claim}\label{xgandy}
The subgroup $[X,G]$ generates $Y = [X, G]Z$ modulo $Z$. 
\end{claim} 
\begin{proof} 
Working with the group $\bar{G}=G/\Ker(\chi)$, using the same argument as in the proof of 
Claim 3.26 of \cite{edith1}, we have that $[\bar{X},\bar{G}]$ generates
$\bar{Y} = [\bar{X}, \bar{G}]\bar{Z}$ modulo $\bar{Z}$. Since $Z={\bf Z}(\chi)$, we have 
that $\Ker(\chi)\leq Z$. Thus 
$\bar{Z}= Z/\Ker(\chi)$ and the claim follows.
\end{proof}

\begin{comment}\label{choosingcoset}
Observe that  $G/H$ is cyclic of order $p$. So we may choose $g \in G$
such that the distinct cosets of $H$ in $G$ are $H$, $Hg$, $Hg^2$, $\ldots$,
$Hg^{p-1}$.
\end{comment}

 Since $\chi= \xi^G$ and $\xi_X=\alpha$, it follows from \ref{xcentralh} that
\begin{equation*}
\chi_X=  \alpha + \alpha^g + \cdots + \alpha^{g^{p-1}}=
\sum_{i=0}^{p-1} \alpha^{g^i} .
\end{equation*}
Similarly, we have that 
\begin{equation*}
\psi_X=  \beta + \beta^g + \cdots + \beta^{g^{p-1}}
=    \sum_{j=0}^{p-1} \beta^{g^j}.
\end {equation*}
Combining the two previous equations we have that
\begin{equation}\label{intermedia}
\chi_X \psi_X=  (\sum_{j=0}^{p-1} \alpha^{g^j}))( \sum_{j=0}^{p-1} \beta^{g^j})= 
\sum_{i=0}^{p-1}\sum_{j=0}^{p-1} \alpha^{g^i} \beta^{g^j}. 
\end{equation}
By \eqref{product} and \eqref{thetaidelta}, we have that
\begin{equation}\label{comparing}
 (\chi\psi)_X = (\sum_{i=1}^n m_i \theta_i)_X=
\sum_{i=1}^n m_i [\sum_{j=0}^{p-1}(\delta_i\alpha\beta)^{g^j}].
\end{equation}

\begin{claim}\label{orbits} Let $g\in G$ be 
 as in \ref{choosingcoset}.
For each $i=0,1,\dots ,p-1$, there exist $j \in \{0,1,\ldots, p-1\}$
and $\delta_{g^i} \in \Lin( X \modu Y)$  such that 
\begin{equation}\label{productitermsdelta}
\alpha \beta^{g^i} = (\alpha\beta)^{g^j} \delta_{g^i}.
\end{equation}
Also $|\{\delta_{g^i} \mid i=0,1,2,\ldots, p-1\}| \leq n$.
\end{claim}
\begin{proof} 
See Proof of Claim 3.30 of \cite{edith1}.
\end{proof}

\begin{claim}\label{choosingright} 
Let $g \in G$ be
 as in \ref{choosingcoset}.
Then there exist three distinct integers 
 $i,j, k \in \{0,1,2,\ldots, p-1 \}$,
and some $\delta \in \Irr(X \modu Y)$, such that 
\begin{eqnarray*}
\alpha \beta^{g^i} & = & (\alpha\beta)^{g^r} \delta, \\ 
\alpha \beta^{g^j} & = &  (\alpha\beta)^{g^s} \delta, \\
\alpha \beta^{g^k} & = &  (\alpha\beta)^{g^t} \delta,
\end{eqnarray*}

\noindent for some $r,s,t \in \{0,1,2,\ldots, p-1 \}$.
\end{claim}
\begin{proof}
See Proof of Claim 3.34 of \cite{edith1}.
\end{proof}

\begin{claim}\label{choosingcosetgood}
We can choose the element  $g$ in \ref{choosingcoset}
such that one of the following 
holds:

(i) There exists some $j=2, \ldots, p-1$ such that
\begin{eqnarray*}
\alpha \beta^{g} & = & (\alpha\beta)^{g^r},\\
\alpha \beta^{g^j} & = &  (\alpha\beta)^{g^s},
\end{eqnarray*}
\noindent for some $r,s \in \{0,1,\ldots, p-1\}$ with $r \neq 1$.

(ii) There exist $j$ and $k$ such that $1<j<k<p$, and 
\begin{eqnarray*}\label{1,j,k}
\alpha \beta^{g} & = & (\alpha\beta)^{g^r}\delta,\\ 
\alpha \beta^{g^j} & = & (\alpha\beta)^{g^s}\delta, \\
\alpha \beta^{g^k} & = &  (\alpha\beta)^{g^t}\delta,
\end{eqnarray*}
\noindent for some  $\delta \in \Irr( X \modu Y)$ and some 
$r,s,t\in \{0,1,\ldots, p-1\}$ with $r\neq 1$.

\end{claim}
\begin{proof} 
See Proof of Claim 3.35 of \cite{edith1}.
\end{proof}

Let $g$ be as in Claim \ref{choosingcosetgood}.
Since $X/Y$ is cyclic of order $p$, we may choose $x \in X$ such 
that $X = Y<x>$. Since $H$ is abelian, we have $[X,H]=1$. 
Suppose that $[x,g^{-1}]\in Z$. Then $x$ centralizes both 
$g^{-1}$ and $H$ modulo $Z$. Hence $xZ \in {\bf Z}(G/Z)$ and so
$[x,G]\leq Z$. Since $Y/Z$ is a chief
section of the $p$-group $G$, we have that $[Y,G]\leq Z$ and so $[<x>Y, G]=[X,G]\leq Z$ 
which is false by Claim \ref{xgandy}. Hence 
$[x, g^{-1}] \in Y\setminus Z$ and so  
\begin{equation}\label{18}
Y=Z<y> \mbox{ is generated over } Z \mbox{ by } y=[x,g^{-1}].
\end{equation}
Since $[Y,G]\leq Z$ we have that $z=[y,g^{-1}]  \in Z$. If
$z=1$, then $G=H<g>$ centralizes $Y = Z<y>$, since $H$ centralizes 
$Y<X$ by \ref{xcentralh}, and $G$ centralizes Z. This
is impossible because $Z={\bf Z}(G)< Y$. Thus 
\begin{equation}\label{19}
z= [y,g^{-1}] \mbox{ is a non-trivial element of }Z.
\end{equation}
By \eqref{18} we have  $y=[x,g^{-1}]=x^{-1} x^{g^{-1}}$. By
\eqref{19} we have $z= [y, g^{-1}]=y^{-1} y^{g^{-1}}$ . Finally
$z^{g^{-1}}= z$ since 
$z\in Z$. Since $X=Z<x,y> \leq H$ is abelian,
it follows that 
\begin{equation}\label{21}
z^{g^{-j}}=z, \  y^{g^{-j}}= yz^j \mbox{ and }
 x^{g^{-j}}=xy^j z^{\binom{j}{2}},
\end{equation}
\noindent for any integer $j= 0, 1, \ldots, p-1$. 
Because $g^{-p} \in H$ centralizes $X$ by \ref{xcentralh}, 
we have
\begin{equation*} 
z^p=1 \mbox{ and } y^p z^{\binom{p}{2}}=1.
\end{equation*}
Observe that the statement is true for $p\leq 3$ since then $\frac{p+1}{2}\leq 2$. Thus we may
assume that $p$ is odd. Hence $p$ divides
$\binom{p}{2}= \frac{p(p-1)}{2}$ and $z^{\binom{p}{2}}=1$. Therefore

\begin{equation}\label{22}
y^p= z^p =1.
\end{equation}

It follows that $y^i$, $z^i$ and $z^{\binom{i}{2}}$ 
depend only on the residue of $i$ modulo 
$p$, for any integer $i\geq 0$.

\begin{comment}\label{betanot1}
Observe that
$\Ker(\xi_Z) \cap \Ker(\tau_Z) \le \Ker(\chi) \cap \Ker(\psi) = 1$ implies that $z$ is not in both $\Ker(\xi_Z)$ and $\Ker(\tau_Z)$. 
Without loss of generality we may assume that $\tau_Z(z) \neq 1$. Since $\beta$ is an
extension of $\tau_Z$, we may assume that  $\beta(z)\neq 1$.
\end{comment}

\begin{claim}\label{nottrivialp}
$\xi_Z\tau_Z(z)$ is primitive $p$-th root of unit.
\end{claim}
\begin{proof}
Suppose that $(\xi_Z\tau_Z)(z)=1$. Then $(\xi_Z\tau_Z)([y,g^{-1}])=1$ and so $(\xi_Z\tau_Z)^g(y)=(\xi_Z\tau_Z)(y)$.
 Since $H$ is abelian, $|G:H|=p$, 
$\theta_i$ lies above $\xi_Z\tau_Z$ for all $i$,  
and $g\in G\setminus H$, it follows that $Y=<y, {\bf Z}(G)>$ is 
contained in ${\bf Z}(\theta_i)$.
That is  contradiction with Claim \ref{faithful}. Thus $(\xi_Z\tau_Z)(z)\neq 1$. Since $z$ is of order $p$ and $\xi_Z\tau_Z$ is a linear character, the claim follows. 
\end{proof}
\begin{claim}\label{contradiction}
 Suppose that
\begin{equation}\label{23}
\alpha \beta^{g} = (\alpha\beta)^{g^r}\delta, 
\end{equation}
\noindent and 
\begin{equation}\label{26}
\alpha \beta^{g^j}  = (\alpha\beta)^{g^s} \delta,
\end{equation}
\noindent  for some 
$j \in \{0, 1, \ldots, p-1\}$, $j \neq 1 $, 
 some $\delta \in \Irr(X \modu Y)$, and some 
$r, s\in \{0,1, \ldots , p-1\}$.
Then           
\begin{equation}\label{contradiction2}
\delta(x)= \beta(z)^{hj(r-1)},
\end{equation}
\noindent where $2h\equiv 1 \modu p$. 
\end{claim}
\begin{proof}
By Claim \ref{nottrivialp} and the same argument as in the 
 proof of Claim 3.40 of \cite{edith1}, the statement follow.
\end{proof}

Suppose that 
Claim \ref{choosingcosetgood} (ii) holds. Then by Claim \ref{contradiction}
 we have that 
 $\delta(x)= \beta(z)^{hj(r-1)}$ and
$\delta(x) =\beta(z)^{hk(r-1)}$. Since $\beta(z)=\tau_Z(z)$ is a primitive $p$-th root of unit
by \ref{betanot1}, 
we have that $ hj(r-1)\equiv hk(r-1) \modu p$.
Since $r\not \equiv 1 \modu p$ and $2h\equiv 1\modu p$, we have that  
$k \equiv j\modu p $, a contradiction. Thus 
Claim \ref{choosingcosetgood}
(i) must hold. 

We apply now Claim \ref{contradiction} with  $\delta=1$. 
Thus $1= \delta(x) =  \beta(z)^{hj(r-1)}$.
Therefore $hj(r-1)\equiv 0 \modu p$. Since
$2h\equiv 1 \modu p$, either $j\equiv 0 \modu p$ or $r-1 \equiv 0 \modu p$. 
Neither is possible. That is our final contradiction and Lemma \ref{basicolemma}
is proved.
\end{proof}
\begin{proof}[Proof of Theorem A]
Since $G$ is a nilpotent group, $G$ is the direct product $G_1 \times G_2$ of its Sylow $p$-subgroup $G_1$ and its Hall $p'$-subgroup $G_2$. We can write then $\chi=\chi_1\times\chi_2$ and $\psi=\psi_1\times\psi_2$ for some characters $\chi_1,\psi_1\in \Irr(G_1)$ and some characters
$\chi_2,\psi_2\in \Irr(G_2)$. Since $\chi(1)=p$, we have that
$\chi_2(1)=1$ and thus $\chi_2\psi_2\in \Irr(G_2)$. 
If $\psi(1)\neq p$, since $\psi(1)$ is a prime number, we have that $\psi_1(1)=1$ and thus 
 $\chi_1\psi_1$ is an irreducible. Therefore $\chi\psi\in \Irr(G)$ and (iv) holds.
 We may assume then that $\psi(1)=p$ and thus $\psi_2(1)=1$. Then $\chi_2\psi_2$ is 
 a linear character and so we may assume that
  $G$ is a $p$-group. 

If $\chi\psi$ has a linear constituent, by Lemma \ref{derived2} and Lemma \ref{linearconstituent}
we have that (i) or (ii) holds. So we may assume that all the irreducible constituents of 
$\chi\psi$ are of degree at least $p$. If $\chi\psi$ has an irreducible constituent of degree 
$p^2$, then $\chi\psi\in \Irr(G)$ and (iv) holds. We may assume then that all the irreducible
constituents of $\chi\psi$ have degree $p$. Since $\chi\psi(1)=p^2$, it follows that
$\eta(\chi\psi)\leq p$. By Lemma \ref{basicolemma} we have that either
$\eta(\chi\psi)=1$ or 
$\eta(\chi\psi)\geq \frac{p+1}{2}$, and so (iii) holds. 
\end{proof}

{\bf Examples.} Fix a prime $p>2$

{\bf (i)}
Let $E$ be an extraspecial group of order $p^3$ and $\phi\in \Irr(E)$ of degree $p$. We can check that  the product
$\phi\overline{\phi}$ is the sum of all the linear characters of $E$.

{\bf (ii)}
In the proof of Proposition 6.1 of \cite{edith2}, an example is constructed of a $p$-group 
$G$ and a character $\chi\in \Irr(G)$ such that $\chi\overline{\chi}$ is the sum of $p$ distinct linear characters and of $p-1$ distinct irreducible characters of degree $p$.

{\bf (iii)}
Given an extraspecial group $E$ of order $p^3$, where $p>2$, and $\phi\in \Irr(E)$ a character of degree $p$, we can check  that 
$\phi\phi$ is a multiple of an irreducible.  
In Proposition 6.1 of \cite{edith1}, an example 
is provided of a $p$-group $G$ and a character $\chi\in \Irr(G)$ such that
$\eta(\chi\chi)=\frac{p+1}{2}$. In \cite{maria}, an example is provided of a $p$-group 
$P$
and two faithful characters $\delta,\epsilon\in \Irr(P)$ of degree $p$ such that 
$\eta(\delta\epsilon)=p-1$. 
     
Let $G$ be the wreath product of a group of order $p$ with itself. Thus $G$ has a normal abelian 
subgroup of $N$ of
order $p^p$ and index $p$. Let $\lambda\in \Lin(N)$. We can check that
$\chi=\lambda^G$ and $\psi=(\lambda^2)^G$ are irreducible characters of degree
$p$ and  $\chi\psi$ is the sum of $p$ distinct irreducible characters of degree 
$p$.

We wonder if there exists a $p$-group $P$ with characters $\chi,\psi\in \Irr(P)$ of
degree $p$ such that 
$\frac{p+1}{2}<\eta(\chi\psi)<p-1$.

{\bf (iv)} Let $Q$ be a $p$-group and $\kappa\in \Irr(Q)$ be a character of degree $p$.
Set $P=Q\times Q$, $\chi=\kappa\times 1_G$ and $\psi=1_G\times \kappa$. Observe that
$\chi$, $\psi$ and $\chi\psi$ are irreducible characters of $P$.

\end{section}
\bibliographystyle{amsplain}

\begin{thebibliography}{9}

\bibitem{edith1}
E. Adan-Bante, Products of characters and finite p-groups,  J. of Algebra, 277 (2004) 236-255.

\bibitem{edith2}
E. Adan-Bante, Products of characters and finite p-groups II, Arch. Math. 82  No 4 (2004), 289-297.


\bibitem{edithderived2}
E. Adan-Bante, Products of characters and derived length II, J. Group Theory 8 (2005), 453-459.

\bibitem{restriction} E. Adan-Bante, Products of characters and restriction of characters, 
preprint.

\bibitem{isaacs} I.M.Isaacs, Character Theory of Finite Groups. New York-San
Francisco--London: Academic Press 1976

\bibitem{maria}
M. Loukaki, A. Moreto, 
On the number of constituents of products of characters , 
Algebra Colloq. 14, no.2, 207-208,(2007).
\end{thebibliography}

\end{document}